\documentclass[submission%
]{dmtcs}

\usepackage[latin1]{inputenc}
\usepackage{subfigure}

\usepackage[round]{natbib}

\usepackage[all]{xy}
\usepackage{amsmath}
\usepackage{makeidx}

\usepackage{tikz}
\usetikzlibrary{decorations.pathreplacing}
\usetikzlibrary{matrix,arrows}

\numberwithin{equation}{section}
\newtheorem{theorem}[subsection]{Theorem}
\newtheorem{definition}[subsection]{Definition}

\newcommand{\Z}{\mathbb{Z}}
\newcommand{\N}{\mathbb{N}}
\newcommand{\Q}{\mathbb{Q}}

\newcommand{\ce}{{\mathcal E}}

\newcommand{\A}{\mathbb{A}}
\newcommand{\E}{\mathbb{E}}

\newcommand{\eps}{\varepsilon}

\newcommand{\sgn}{\mbox{sign}}
\newcommand{\ul}[1]{\bold{#1}}

\newcommand{\Frac}{\mbox{Frac}}
\newcommand{\Ad}{\mbox{Ad}}
\newcommand{\Id}{\mbox{Id}}

\newcommand{\iso}{\stackrel{_\sim}{\rightarrow}}

\newcommand{\bt}{\bullet}

\author{Bernhard Keller\addressmark{1}\thanks{Email: \email{keller@math.jussieu.fr}}}\title{Quiver mutation and combinatorial DT-invariants}
\address{\addressmark{1}Universit\'e Paris Diderot -- Paris~7, Institut
Universitaire de France,  Institut de Math\'ematiques de Jussieu -- Paris Rive Gauche, 
UMR 7586 du CNRS, Paris, France}
\keywords{Cluster algebra, quiver mutation, Donaldson-Thomas invariants}
\received{2013-05-26}
\begin{document}
\maketitle
\begin{abstract}
A quiver is an oriented graph. Quiver mutation is an elementary
operation on quivers. It appeared in physics in Seiberg duality in
the nineties and in mathematics in the definition of cluster algebras
by Fomin-Zelevinsky in 2002. We show how, for large
classes of quivers $Q$, using quiver mutation and quantum dilogarithms,
one can construct the combinatorial DT-invariant, a formal power
series intrinsically associated with $Q$. When defined, it coincides
with the `total' Donaldson-Thomas invariant of $Q$ (with a generic
potential) provided by algebraic geometry (work of Joyce, 
Kontsevich-Soibelman, Szendroi and many others). We
illustrate combinatorial DT-invariants on many examples 
and point out their links to quantum cluster algebras and 
to (infinite) generalized associahedra.

Un carquois est un graphe orient\'e. La mutation des carquois est
une op\'eration \'el\'ementaire sur les carquois. Elle est apparue en
physique dans la dualit\'e de Seiberg dans les ann\'ees 90 et en
math\'ematiques dans la d\'efinition des alg\`ebres amass\'ees par
Fomin--Zelevinsky en 2002. Nous montrons comment, pour de
grandes classes de carquois $Q$, \`a l'aide de la mutation des
carquois et des dilogarithmes quantiques, on peut construire
l'invariant DT combinatoire, une s\'erie formelle associ\'ee 
intrins\`equement \`a $Q$.  Quand cet invariant est d\'efini,
il est \'egal \`a l'invariant de Donaldson--Thomas `total' 
associ\'e \`a $Q$ (avec un potentiel g\'en\'erique) qui est
fourni par la g\'eom\'etrie alg\'ebrique (travaux de Joyce,
Kontsevich--Soibelman, Szendroi et beaucoup d'autres).
Nous illustrons les invariants DT combinatoires sur beaucoup
d'exemple et \'evoquons leurs liens avec les alg\`ebres
amass\'ees quantiques et des associa\`edres g\'en\'eralis\'es
(infinis).
\end{abstract}

\section{Introduction}
\label{s:introduction}

A quiver is an oriented graph. Quiver mutation is an elementary
operation on quivers. It appeared in physics already in the nineties
in Seiberg duality, cf.  \cite{Seiberg95}. In mathematics, quiver
mutation was introduced by  \cite{FominZelevinsky02} as the basic
combinatorial ingredient of their definition of cluster algebras.
Thus, quiver mutation is linked to the large array of subjects
where cluster algebras have subsequently turned out to be
relevant, cf. for example the cluster algebras portal maintained by
\cite{Fomin07} and the survey articles by \cite{Fomin10},
\cite{Leclerc10}, \cite{Reiten10a}, \cite{Williams12}. Among these links, 
the one to representation
theory and algebraic geometry has been particularly
fruitful. It has allowed to `categorify' cluster algebras and thereby
to prove conjectures about them which seem beyond
the scope of the purely combinatorial methods, cf.
for example the articles of 
\cite{DerksenWeymanZelevinsky10}
(which ultimately builds on \cite{CalderoChapoton06}), 
\cite{GeissLeclercSchroeer11b},
\cite{Plamondon11a},
\cite{CerulliKellerLabardiniPlamondon12},  \ldots\ .
  
The constructions and results we present in this talk are another
manifestation of this fruitful interaction. They are inspired by
the theory of Donaldson--Thomas invariants as it has been
developped by Bridgeland, 
\cite{JoyceSong09},
\cite{KontsevichSoibelman08, KontsevichSoibelman10},
\cite{Nagao10},
\cite{Reineke11},
\cite{Szendroi08} and many others.
In this theory, one assigns Donaldson--Thomas invariants
to three-dimensional, possibly non commutative, Calabi--Yau varieties.
These invariants exist in many different versions. Here we use
the `total' Donaldson--Thomas invariant, which is a certain
power series in (slighly) non commutative variables. One
important construction of non commutative $3$-Calabi--Yau
varieties takes as its input a quiver (with a generic potential). Thus,
there is a Donaldson--Thomas invariant associated with 
`each' quiver (some technical problems remain to be
solved for a completely general definition). It
turns out that for a suprisingly large class of quivers, it is
possible to construct this invariant in a combinatorial way
using products of quantum dilogarithm series associated
with so-called reddening sequences of quiver mutations.
This construction yields the definition of the combinatorial
DT-invariant, which is the main point of this talk 
(section~\ref{s:combinatorial-DT-invariants}). 
It is an important fact that a given quiver may admit many distinct
reddening sequences. Each of them yields a product
decomposition for the combinatorial DT-invariant and
in this way, one obtains many interesting quantum dilogarithm
identities (section~\ref{s:examples}).

Let us emphasize that the use of (products of) quantum dilogarithm
series in the study of (quantum) cluster algebras goes back
to the insight of \cite{FockGoncharov09, FockGoncharov09a}. They
also pioneered their application in the study of (quantum) dilogarithm
identities, which was subsequently developped by many authors.
We refer to \cite{Nakanishi12} for a survey. Geometric as well as
combinatorial constructions of DT-invariants also appear in physics,
cf. for example 
\cite{CecottiNeitzkeVafa10},
\cite{CecottiEtAl11, CecottiEtAl11a}, 
\cite{CecottiCordovaVafa11},
\cite{CecottiVafa11},
\cite{GaiottoMooreNeitzke09, GaiottoMooreNeitzke10,  GaiottoMooreNeitzke10a},
\cite{Xie12},
\ldots .

\section{Quiver mutation}
\label{s:quiver-mutation}

A {\em quiver}\index{quiver} is an oriented graph, i.e. a
quadruple $Q=(Q_0, Q_1, s, t)$ formed by a set of vertices $Q_0$, a set
of arrows $Q_1$ and two maps $s$ and $t$ from $Q_1$ to $Q_0$ which send
an arrow $\alpha$ respectively to its source $s(\alpha)$ and its
target $t(\alpha)$. In practice, a quiver is given by a picture
as in the following example
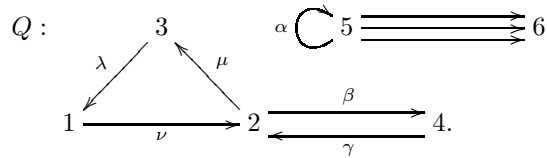
\[ Q:
\xymatrix{ & 3 \ar[ld]_\lambda & & 5 \ar@(dl,ul)[]^\alpha \ar@<1ex>[rr] \ar[rr] \ar@<-1ex>[rr] & & 6 \\
  1 \ar[rr]_\nu & & 2 \ar@<1ex>[rr]^\beta \ar[ul]_\mu & & 4.
  \ar@<1ex>[ll]^\gamma }
\]
An arrow $\alpha$ whose source and target coincide is a {\em loop}\index{loop}; a
{\em $2$-cycle}\index{$2$-cycle} is a pair of distinct arrows $\beta$ and $\gamma$ such that
$s(\beta)=t(\gamma)$ and $t(\beta)=s(\gamma)$. Similarly, one defines
{\em $n$-cycles} for any positive integer $n$. A vertex $i$ of a quiver
is a {\em source}\index{source of a quiver} (respectively a {\em sink}\index{sink of a quiver}) if there is no arrow
with target $i$ (respectively with source $i$).
A {\em Dynkin quiver} is a quiver whose underlying graph is a
Dynkin diagram of type $A_n$, $n\geq 1$,  $D_n$, $n\geq 4$, 
or $E_6$, $E_7$, $E_8$.

By convention, in the sequel, by a quiver we always mean a finite
quiver without loops nor $2$-cycles whose set of vertices is the
set of integers from $1$ to $n$ for some $n\geq 1$. Up to an
isomorphism fixing the vertices, such a quiver $Q$ is given by
the {\em skew-symmetric matrix $B=B_Q$} whose coefficient $b_{ij}$ is
the difference between the number of arrows from $i$ to $j$ and
the number of arrows from $j$ to $i$ for all $1\leq i,j\leq n$.
Conversely, each skew-symmetric matrix $B$ with integer coefficients
comes from a quiver.

Let $Q$ be a quiver and $k$ a vertex of $Q$. The {\em mutation $\mu_k(Q)$}\index{mutation!of a quiver}
is the quiver obtained from $Q$ as follows:
\begin{itemize}
\item[1)] for each subquiver $\xymatrix{i \ar[r]^\beta & k \ar[r]^\alpha & j}$,
we add a new arrow $[\alpha\beta]: i \to j$;
\item[2)] we reverse all arrows with source or target $k$;
\item[3)] we remove the arrows in a maximal set of pairwise disjoint $2$-cycles.
\end{itemize}
For example, if $k$ is a source or a sink of $Q$, then the
mutation at $k$ simply reverses all the arrows incident with $k$. In general,
if $B$ is the skew-symmetric matrix associated with $Q$ and $B'$ the one
associated with $\mu_k(Q)$, we have
\begin{equation} \label{eq:matrix-mutation}
b'_{ij} = \left\{ \begin{array}{ll}
-b_{ij} & \mbox{if $i=k$ or $j=k$~;} \\
b_{ij}+\sgn(b_{ik})\max(0, b_{ik} b_{kj}) & \mbox{else.}
\end{array} \right.
\end{equation}
This is the {\em matrix mutation rule}\index{mutation!of a matrix} for skew-symmetric
(more generally: skew-symmetri\-zable) matrices introduced by
\cite{FominZelevinsky02}, cf. also \cite{FominZelevinsky07}.

One checks easily that $\mu_k$ is an involution. For example,
the quivers
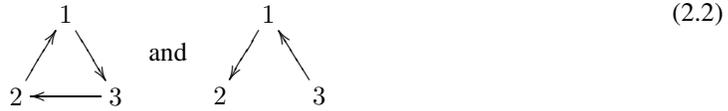
\begin{equation} \label{quiver1}
\begin{xy} 0;<0.2pt,0pt>:<0pt,-0.2pt>::
(94,0) *+{1} ="0",
(0,156) *+{2} ="1",
(188,156) *+{3} ="2",
"1", {\ar"0"},
"0", {\ar"2"},
"2", {\ar"1"},
\end{xy}
\begin{minipage}{1cm}
\vspace*{1cm}
\begin{center} and \end{center}
\end{minipage}
\begin{xy} 0;<0.2pt,0pt>:<0pt,-0.2pt>::
(92,0) *+{1} ="0",
(0,155) *+{2} ="1",
(188,155) *+{3} ="2",
"0", {\ar"1"},
"2", {\ar"0"},
\end{xy}
\end{equation}
are linked by a mutation at the vertex $1$. Notice that these
quivers are drastically different: The first one is a cycle,
the second one the Hasse diagram of a linearly ordered set.

Two quivers are {\em mutation equivalent}\index{mutation equivalent} if they are linked
by a finite sequence of mutations. For example, it is an easy
exercise to check that any two orientations of a tree are
mutation equivalent. Using the quiver mutation applet
by \cite{KellerQuiverMutationApplet}
or the Sage package by \cite{MusikerStump11} one can check
that the following three quivers are mutation equivalent
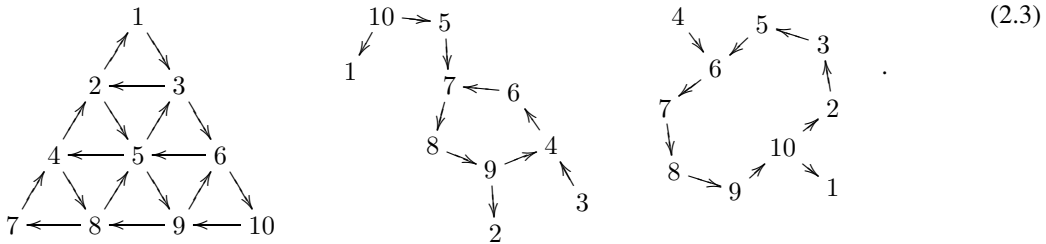
\begin{equation} \label{quiver3}
\begin{xy} 0;<0.6pt,0pt>:<0pt,-0.6pt>::
(79,0) *+{1} ="0",
(52,44) *+{2} ="1",
(105,44) *+{3} ="2",
(26,88) *+{4} ="3",
(79,88) *+{5} ="4",
(131,88) *+{6} ="5",
(0,132) *+{7} ="6",
(52,132) *+{8} ="7",
(105,132) *+{9} ="8",
(157,132) *+{10} ="9",
"1", {\ar"0"},
"0", {\ar"2"},
"2", {\ar"1"},
"3", {\ar"1"},
"1", {\ar"4"},
"4", {\ar"2"},
"2", {\ar"5"},
"4", {\ar"3"},
"6", {\ar"3"},
"3", {\ar"7"},
"5", {\ar"4"},
"7", {\ar"4"},
"4", {\ar"8"},
"8", {\ar"5"},
"5", {\ar"9"},
"7", {\ar"6"},
"8", {\ar"7"},
"9", {\ar"8"},
\end{xy}
\quad\quad
\begin{xy} 0;<0.3pt,0pt>:<0pt,-0.3pt>::
(0,70) *+{1} ="0",
(183,274) *+{2} ="1",
(293,235) *+{3} ="2",
(253,164) *+{4} ="3",
(119,8) *+{5} ="4",
(206,96) *+{6} ="5",
(125,88) *+{7} ="6",
(104,164) *+{8} ="7",
(177,194) *+{9} ="8",
(39,0) *+{10} ="9",
"9", {\ar"0"},
"8", {\ar"1"},
"2", {\ar"3"},
"3", {\ar"5"},
"8", {\ar"3"},
"4", {\ar"6"},
"9", {\ar"4"},
"5", {\ar"6"},
"6", {\ar"7"},
"7", {\ar"8"},
\end{xy}
\quad\quad
\begin{xy} 0;<0.3pt,0pt>:<0pt,-0.3pt>::
(212,217) *+{1} ="0",
(212,116) *+{2} ="1",
(200,36) *+{3} ="2",
(17,0) *+{4} ="3",
(123,11) *+{5} ="4",
(64,66) *+{6} ="5",
(0,116) *+{7} ="6",
(12,196) *+{8} ="7",
(89,221) *+{9} ="8",
(149,166) *+{10} ="9",
"9", {\ar"0"},
"1", {\ar"2"},
"9", {\ar"1"},
"2", {\ar"4"},
"3", {\ar"5"},
"4", {\ar"5"},
"5", {\ar"6"},
"6", {\ar"7"},
"7", {\ar"8"},
"8", {\ar"9"},
\end{xy}
\begin{minipage}{1cm}
\vspace*{1.5cm}
\begin{center} . \end{center}
\end{minipage}
\end{equation}
The common {\em mutation class}\index{mutation class} of these quivers contains 5739 quivers
(up to isomorphism). The mutation class of `most' quivers is infinite.
The classification of the quivers having a finite mutation class
was achieved by by 
\cite{FeliksonShapiroTumarkin12} (and by \cite{FeliksonShapiroTumarkin12a}
in the skew-symmetric case):
in addition to the quivers associated
with triangulations of surfaces (with boundary and marked points,
cf. \cite{FominShapiroThurston08}),
the list contains $11$ exceptional quivers, the largest of which is
in the mutation class of the quivers~(\ref{quiver3}).

\section{Green quiver mutation} 

Let $Q$ be a quiver without loops nor $2$-cycles.
The \emph{framed quiver} $\tilde{Q}$ is
obtained from $Q$ by adding, for each vertex $i$, a new
vertex $i'$ and a new arrow $i \to i'$.
Here is an example:
\[
Q: \xymatrix@R=0.3cm@C=0.3cm{1 \ar[r] & 2} \quad\quad\quad
\tilde{Q} : \xymatrix@R=0.3cm@C=0.3cm{ 1 \ar[r] \ar[d] & 2 \ar[d] \\ 1' & 2'}\quad .
\]
The vertices $i'$ are called {\em frozen vertices} because we
never mutate at them. Now suppose that we have transformed 
$\tilde{Q}$ into $\tilde{Q}\,'$ by a finite sequence of mutations 
(at non frozen vertices). A vertex $i$ of $Q$ is \emph{green} in 
$\tilde{Q}\,'$ if there are no arrows ${j\,'} \to i$ in $\tilde{Q}\,'$. 
Otherwise, it is \emph{red}.
A sequence $\ul{i}=(i_1,\ldots, i_N)$ is
\emph{green} if for each $1\leq t\leq N$, the vertex $i_t$
is green in the partially mutated quiver $\tilde{Q}(\ul{i},t)$
defined by
\[
\tilde{Q}(\ul{i},t)= \mu_{i_{t-1}} \ldots \mu_{i_2}\mu_{i_1} (\tilde{Q}),
\]
where for $t=1$, we have the empty mutation sequence and
obtain the initial quiver $\tilde{Q}$. It is {\em maximal green}
if it is green and all the vertices of the final quiver $\mu_{\ul{i}}(\tilde{Q})$
are red (so that indeed, the sequence $\ul{i}$ cannot be extended to any strictly
longer green sequence). 

\begin{figure}
\[
\def\g#1{\save [].[dr]!C *++\frm{}="g#1" \restore}
\xymatrix{
 & &  *+[o][F-]{1}\g1 \ar[r] \ar[d] & *+[o][F-]{2} \ar[d] & & *+[o][F-]{1}\g2  \ar[d] \ar[dr] & 2 \ar[l] & &  \\
 & &   1'                & 2' & &    1'               & 2' \ar[u]\\
\g3 1 & *+[o][F-]{2} \ar[l] \ar[d]  &  &  &   &    &   & \g4 1 \ar[r] & *+[o][F-]{2} \ar[dl] \\
1'\ar[u] & 2'            &  &  &   &    &   &  1' \ar[u] & 2' \ar[ul] \\
 & & \g5 1 \ar[r] & 2         &  &  \g6 1 & 2 \ar[l] \\
 & &    1' \ar[u] & 2' \ar[u] &  &   1' \ar[ur] & 2' \ar[ul]
 \ar@{->}^{\mu_2} "g1"; "g2"
 \ar@<-20pt>@{->}_{\mu_1} "g1"; "g3"
 \ar@<20pt>@{->}^{\mu_1} "g2"; "g4"
 \ar@<-20pt>@{->}_{\mu_2} "g3"; "g5"
 \ar@<20pt>@{->}^{\mu_2} "g4"; "g6"
 \ar@{-}^{\mbox{\small frozen}}_{\mbox{\small isom}} "g5"; "g6"
 \ar_{\mu_{12}} "g1"; "g5"
 \ar^{\mu_{212}} "g1"; "g6"
 }
\]
\caption{The two maximal green sequences for $A_2$}
\label{fig:max-green-sequences}
\end{figure}
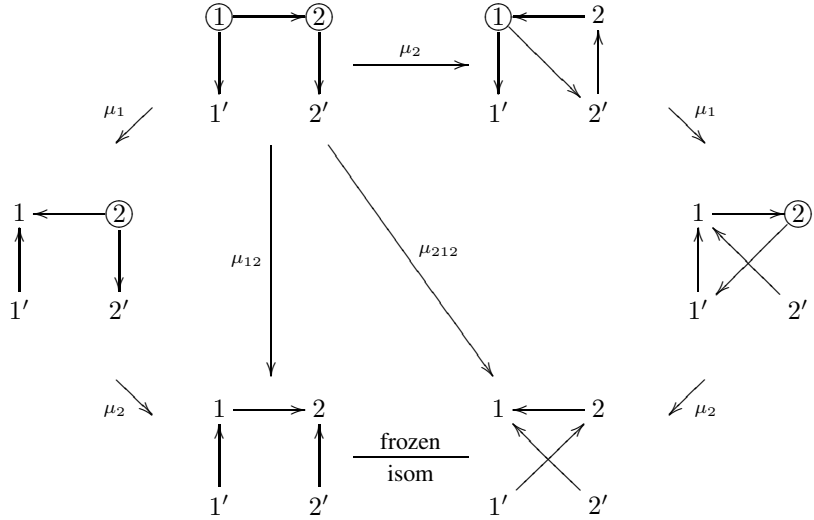

In Figure~\ref{fig:max-green-sequences},
we have encircled the green vertices. We see that in this example,
we have two maximal green sequences $12$ and $212$ and that the final
quivers associated with these two sequences are isomorphic
by an isomorphism which fixes the frozen vertices. We
call such an isomorphism a {\em frozen isomorphism}.
A sequence $\ul{i}=(i_1,\ldots, i_N)$
is \emph{reddening} if all vertices of the final quiver
$
\mu_{\ul{i}}(\tilde{Q}) = \tilde{Q}(\ul{i},N)
$
are red. Of course, maximal green sequences are reddening.
In the example of
Figure~\ref{fig:max-green-sequences}, the sequence
$
1,2,1,2,1,2,1
$
is reddening but not green. 

\begin{theorem} \label{thm:final-quivers}
If $\ul{i}$ and $\ul{i'}$ are reddening sequences,
there is a frozen isomorphism between the final quivers
\[
\mu_{\ul{i}}(\tilde{Q}) \iso \mu_{\ul{i}'}(\tilde{Q}).
\]
\end{theorem}

The statement of the theorem is purely combinatorial but
the known proofs (cf. section~7 of \cite{Keller12} and the references
given there) are based on representation theory and geometry.
For an {\em arbitrary} sequence $\ul{i}$ of non frozen vertices,
we define the $c$-matrix $C(\ul{i})$ as the $n\times n$-matrix
occuring in the right upper corner of the skew-symmetric matrix
associated with the final quiver $\mu_{\ul{i}}(\tilde{Q})$, so that
we have
\[
B_{\mu_{\ul{i}}(\tilde{Q})} = 
\left[\begin{array}{cc} 
* &  C(\ul{i}) \\
* & * 
\end{array} 
\right].
\]
Thus, the $(i,j)$-coefficient of the matrix $C(\ul{i})$ is
the difference between the number of arrows $i \to j'$
and $j' \to i$. The {\em c-vectors} associated with the
sequence $\ul{i}$ are by definition the rows of
the matrix $C(\ul{i})$. The following statement is known
as the {\em sign coherence of c-vectors}.

\begin{theorem}[\cite{DerksenWeymanZelevinsky10}]
\label{thm:sign-coherence}
Each $c$-vector lies in $\N^n$ or $(-\N)^n$.
\end{theorem}

Again, the known proof of this combinatorial statement uses representation 
theory and geometry.
The {\em oriented exchange graph} of the quiver $Q$ is defined to be
the quiver $\ce_Q$ whose vertices are the frozen isomorphism
classes of the quivers $\mu_{\ul{i}}(\tilde{Q})$, where $\ul{i}$ is
an arbitrary sequence of vertices of $Q$, and
where we have an arrow
\[
\tilde{Q}' \to \mu_j(\tilde{Q}')
\]
whenever $j$ is a green vertex of $\tilde{Q}'$. 
For example, if $Q$ is the quiver $1 \to 2$, then we see from
Figure~\ref{fig:max-green-sequences} that the oriented exchange
graph is the oriented pentagon
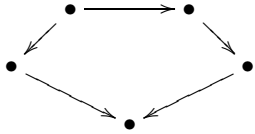
\[
\xymatrix@C=0.4cm@R=0.4cm{  & \bullet \ar[rr] \ar[ld] & & \bullet \ar[rd] & \\
\bullet \ar[rrd] & & & & \bullet \ar[lld] \\
 &  &  \bullet &  & }
\]
By Theorem~\ref{thm:final-quivers}, the quiver $\ce_Q$ has at
most one sink. One can also show that it always has a unique source.
An arbitrary sequence $\ul{i}$ of non frozen vertices corresponds
to a {\em walk} in $\ce_Q$ (a sequence of arrows and formal
inverses of arrows) and a reddening sequence to a {\em path}
(a formal composition of arrows) 
from the source to the sink. If $Q$ is mutation equivalent to a quiver
whose underlying graph is a Dynkin diagram of type $A_n$ 
(resp. $\Delta$), then $\ce_Q$ is
an orientation of the $1$-skeleton of the $n$th Stasheff associahedron
(resp. of the generalized associahedron of type $\Delta$,
cf. \cite{ChapotonFominZelevinsky02}). 

One can show that $\ce_Q$ is the Hasse graph of a poset
(a subposet of the set of torsion subcategories, cf. section~7.7
of \cite{Keller12}). If $Q$ is the linear orientation of $A_n$, this poset
is the Tamari lattice, by \cite{Thomas11}, and for other quivers, one obtains
Cambrian lattices in the sense of \cite{Reading06}.
For certain classes of quivers, this poset is
studied from the viewpoint of representation theory for example by
\cite{AdachiIyamaReiten12}, \cite{BruestleDupontPerotin12},
\cite{KingQiu11}, \cite{KoenigYang12} and
\cite{Ladkani07}. 
Not all quivers admit reddening sequences. For example the
quiver
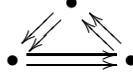
\[
\xymatrix@R=0.4cm@C=0.4cm{ & \bullet \ar@<-2pt>[dl] \ar@<2pt>[dl] & \\
\bullet \ar@<2pt>[rr] \ar@<-2pt>[rr] &  & \bullet \ar@<2pt>[ul] \ar@<-2pt>[ul] }
\]
does not admit a reddening sequence. On the other hand, reddening
sequences do exist
for large classes of quivers, cf. section~\ref{s:examples} below. In particular, each acyclic
quiver (=quiver without oriented cycles) admits a maximal green sequence
corresponding to an increasing enumeration of the vertices for the
order defined by the existence of a path. 

\section{Combinatorial DT-invariants}
\label{s:combinatorial-DT-invariants}

Our aim is to associate an intrinsic formal power series $\E_Q$
with each quiver $Q$ admitting a reddening sequence. 
For quiver with a unique vertex and no arrows, this series will be
the \emph{quantum dilogarithm series} 
\begin{align*}
\E(y)  &\ = 1 + \frac{q^{1/2}}{q-1}\cdot y + \cdots  +
\frac{q^{n^{\,2}/2} y^n}{(q^n-1)(q^n-q) \cdots (q^n-q^{n-1})} + \cdots \\
       &\ \in \Q(q^{1/2})[[y]] ,
\end{align*}
where $q^{1/2}$ is an indeterminate whose square is denoted by $q$ and
$y$ is an indeterminate. This series is a classical object with many
remarkable properties, cf. for example \cite{Zagier91}. 
We will focus on one of them, namely the {\em pentagon identity}: If
$y_1$ and $y_2$ are two indeterminates which $q$-commute,
i.e. $y_1 y_2 = q y_2 y_1$, then we have
\begin{equation} \label{eq:pentagon}
\E(y_1) \E(y_2) = \E(y_2) \E(q^{-1/2} y_1 y_2) \E(y_1).
\end{equation}
It is due to \cite{FaddeevVolkov93} and \cite{FaddeevKashaev94}; 
a recent account can be found in \cite{Volkov12}.
Notice a striking structural similarity between this identity and
the diagram in Figure~\ref{fig:max-green-sequences}: The two
factors $\E(y_1) \E(y_2)$ on the left correspond to the two
mutations in the path on the left, the three factors $\E(y_2) \E(q^{-1/2} y_1 y_2)\E(y_1)$
on the right correspond to the three mutations in the path on the
right and the equality corresponds to the frozen isomorphism between
the final quivers. The common value of the two products will be
defined as the combinatorial DT-invariant $\E_Q$ associated
with the quiver $Q: 1 \to 2$. 

Let $Q$ be a quiver.
For {\em any} sequence $\ul{i}=(i_1, \ldots, i_N)$ of vertices of $Q$,
we will define a {\em product $\E_{Q,\ul{i}}$} of quantum dilogarithm
series. For this, let $\lambda_Q:\Z^n \times \Z^n \to \Z$ be 
the bilinear antisymmetric form associated with the matrix $B_Q$.
Define the {\em complete quantum affine space} as the algebra
\[
\hat{\A}_Q = \Q(q^{1/2})\langle\langle y^\alpha, \alpha\in\N^n \;|\;
y^\alpha y^\beta = q^{1/2 \,\lambda(\alpha,\beta)} y^{\alpha+\beta} \rangle\rangle.
\]
This is a slightly non commutative deformation of an ordinary commutative power
series algebra in $n$ indeterminates. Denote by $\A_Q$ the non completed version.
We define the \emph{product}
\[
\E_{Q,\ul{i}} = \E(y^{\eps_1 \beta_1})^{\eps_1} \cdots \E(y^{\eps_N \beta_N})^{\eps_N} ,
\]
where the product is taken in $\hat{\A}_Q$, the vector $\beta_t$ is
the $t$-th column of the $c$-matrix $C(i_1, \ldots, i_{t-1})$ and
$\eps_t$ is the common sign of the entries of this column
(Theorem~\ref{thm:sign-coherence}), $1 \leq t \leq N$.

\begin{theorem} \label{thm:identity} If $\ul{i}$ and $\ul{i}'$ are
two sequences of vertices of $Q$ such that there is a frozen isomorphism
between $\mu_{\ul{i}}(\tilde{Q})$ and $\mu_{\ul{i}'}(\tilde{Q})$, then we
have the equality
$
\E_{Q,\ul{i}} = \E_{Q,\ul{i}'}.
$
\end{theorem}

The theorem is proved in section~7.11 of \cite{Keller12}, cf. also
\cite{Nagao11}.
In particular, if $\ul{i}$ and $\ul{i}'$ are two reddening sequences, then
by Theorem~\ref{thm:identity}, the above equality holds. More generally,
the theorem shows that to each vertex of the oriented exchange graph
$\ce_Q$, a canonical power series in $\hat{\A}_Q$ is associated.
The one associated to the unique sink (if it exists) is $\E_Q$ and
any reddening sequence gives a product expansion for $\E_Q$.

\begin{definition} If $Q$ admits a reddening sequence $\ul{i}$, the
{\em combinatorial DT-invariant of $Q$} is defined as
\[
\E_Q = \E_{Q,\ul{i}} \in \hat{\A}_Q.
\]
The {\em adjoint combinatorial DT-invariant of $Q$} is
$
DT_Q = \Sigma\circ \Ad(\E_Q) : \Frac(\A_Q) \to \Frac(\A_Q),
$
where $\Frac(\A_Q)$ is the non commutative field of fractions
of $\A_Q$ (cf. the Appendix of \cite{BerensteinZelevinsky05}) and $\Sigma$ its automorphism determined by
$\Sigma(y^\alpha) = y^{-\alpha}$ for all $\alpha\in\N^n$. 
\end{definition}

For the agreement with the geometrically defined DT-invariant, we refer to
section~7 of \cite{Keller12}. In physics, an equivalent procedure has been
discovered independently, cf. \cite{Xie12} and the references given there.
It is easy to check that for $Q: 1 \to 2$, the above definition 
yields the left and right hand sides of the pentagon identity~(\ref{eq:pentagon})
associated with the two maximal green sequences of 
Figure~\ref{fig:max-green-sequences} so that indeed, the
combinatorial DT-invariant $\E_Q$ equals these two products.
One can show that in this case, the adjoint combinatorial DT-invariant
satisfies
$
(DT_Q)^5 = \Id.
$
Below, we will explore some remarkable generalizations of this
example. 

\section{Examples}
\label{s:examples}

\subsection{Dynkin quivers} Let $Q$ be an {\em alternating} Dynkin
quiver, i.e. a simply laced Dynkin diagram endowed with an orientation 
such that each vertex is a source or a sink, for example
\[
Q = \vec{A}_5 : \xymatrix{\bt & \circ \ar[l] \ar[r] & \bt & \circ \ar[l]\ar[r] & \bt} .
\]
Let $i_+$ be the sequence of all sources $\circ$ and $i_-$ the
sequence of all sinks $\bt$ (in any order). One can show that in
this case, the sequence
$
\ul{i} = i_+ i_-
$
is maximal green and so is
\[
\ul{i}\,' = \underbrace{i_- i_+ i_- \ldots}_{h\mbox{\scriptsize \ factors}} \;\;\; ,
\]
where $h$ is the Coxeter number of the underlying graph of $Q$.
Thus, we have $\E(\ul{i}) = \E(\ul{i}\,')$ and $\E_Q$ is the common value.
The identity $\E(\ul{i}) = \E(\ul{i}\,')$ is due to \cite{Reineke10}. Using
the geometry of the generalized associahedra of \cite{ChapotonFominZelevinsky02}, 
one can show that it is a consequence of the pentagon identity, 
cf. \cite{Qiu11} for another approach. For the adjoint combinatorial DT-invariant, we have
$
DT_Q^{h+2} = \Id.
$
This is closely related to the original form of the periodicity
conjecture of \cite{Zamolodchikov91} proved by \cite{FominZelevinsky03b}.
We refer to \cite{BruestleDupontPerotin12} for the study of
maximal green sequences for more general acyclic quivers.

\subsection{Square products of Dynkin quivers}
\begin{figure}
\[
\xymatrix@C=0.1cm@R=0.5cm{ & & \bt \ar[rrr] &  &  & \circ \ar[ldd]
& & &
\bt \ar[lll] \ar[rrr] & & & \circ \ar[ldd] \\
     \bt \ar[rrr]|!{[dr];[urr]}\hole & &  & \circ \ar[rd] & & &
\bt \ar[lll]|!{[lu];[dll]}\hole \ar[rrr]|!{[rd];[rru]}\hole & & & \circ \ar[rd] & &  \\
           & \circ \ar[lu] \ar[ruu] \ar[d] & & & \bt \ar[lll] \ar[rrr] & & &
\circ \ar[ruu] \ar[lu] \ar[d] & & & \bt \ar[lll] & \\
           & \bt \ar[rrr] & & & \circ \ar[d] \ar[u] & & &
\bt \ar[lll] \ar[rrr] & & & \circ \ar[u] \ar[d] & \\
           & \circ \ar[u] & & & \bt \ar[lll] \ar[rrr] &  & &
\circ \ar[u] & & & \bt \ar[lll] & }
\]
\caption{The quiver $\vec{A}_4 \square \vec{D}_5$}
\label{fig:a4-prod-d5}
\end{figure}
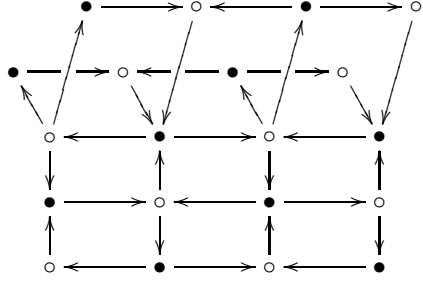

Let $Q_1$ and $Q_2$ be alternating Dynkin quivers and
$Q= Q_1 \square Q_2$ their square product as defined in
section~8 of \cite{Keller10b}. For example,  for
suitable orientations of $A_4$ and $D_5$, the square product
is depicted in Figure~\ref{fig:a4-prod-d5}.
There are no longer sources or sinks in the square product. However, we can consider
the sequence $i_+$ of all even vertices $\circ$ (corresponding to 
a pair of sources or a pair of sinks) and 
the sequence $i_-$ of all odd vertices $\bt$ (corresponding to a mixed pair).
Let 
\begin{align*}
\ul{i}\ \  &= i_+ i_- i_+ \ldots \quad \mbox{with $h$ factors,} \\
\ul{i}\,' &= i_- i_+ i_- \ldots \quad \mbox{with $h\,'$ factors,}
\end{align*}
where $h$ and $h'$ are the Coxeter numbers of the Dynkin diagrams
underlying the two quivers.
One can check that both of these sequences are maximal green.
In particular, the combinatorial DT-invariant is well-defined and
we have
$
\E_Q = \E(\ul{i}) = \E(\ul{i}\,').
$
It is an open question whether these identities are consequences of the
pentagon identity.
One can show (cf. section~5.7 of \cite{Keller11c} or section~8.3.2
of \cite{CecottiNeitzkeVafa10}) that in this case, 
the adjoint combinatorial DT-invariant satisfies
\[
(DT_Q)^m = \Id, \quad \mbox{ where} \quad m = \frac{2(h+h')}{\gcd(h,h')}.
\]

\subsection{Quivers from reduced expressions in Coxeter groups}
If $R$ is an acyclic quiver and $\tilde{w}$ a reduced expression for
an element of the Coxeter group associated with the underlying
graph of $R$, there is a canonical quiver $Q$ associated with the
pair $(R,\tilde{w})$, cf. \cite{BerensteinFominZelevinsky05}.
For example, if $Q$ is $A_4$ with the linear orientation and
$\tilde{w}$ a suitable expression for the longest element, one
obtains the first quiver in (\ref{quiver3}).
As shown by \cite{GeissLeclercSchroeer11b}, such quivers always
admit maximal green sequences. In the $A_4$-example, such a
sequence is given by 
\[
7,8,9,10,4,5,6,2,3,1,7,8,9,2,3,1,7,8,4,7.
\]
Thus, the combinatorial DT-invariant is well-defined. In the above
example (and for all members of the `triangular' family it belongs to), the adjoint combinatorial
DT-invariant satisfies
$
(DT_Q)^6 = \Id.
$
It is an open question for which pairs $(R,\tilde{w})$ the invariant $DT_Q$
is of finite order.

\subsection{Another product construction}
The following quiver is obtained as the triangle product
(cf. section~8 of \cite{Keller10b}) of a quiver of type $A_3$ with
the quiver appearing in the first column (the arrows marked by $2$
are double arrows).
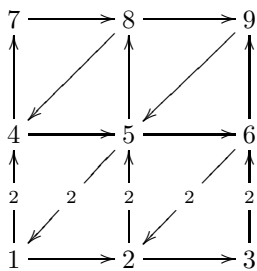
\[
\begin{xy} 0;<0.4pt,0pt>:<0pt,-0.4pt>:: 
(0,227) *+{1} ="0",
(109,227) *+{2} ="1",
(223,227) *+{3} ="2",
(0,109) *+{4} ="3",
(109,109) *+{5} ="4",
(223,109) *+{6} ="5",
(0,0) *+{7} ="6",
(109,0) *+{8} ="7",
(223,0) *+{9} ="8",
"0", {\ar"1"},
"0", {\ar|*+{\scriptstyle 2}"3"},
"4", {\ar|*+{\scriptstyle 2}"0"},
"1", {\ar"2"},
"1", {\ar|*+{\scriptstyle 2}"4"},
"5", {\ar|*+{\scriptstyle 2}"1"},
"2", {\ar|*+{\scriptstyle 2}"5"},
"3", {\ar"4"},
"3", {\ar"6"},
"7", {\ar"3"},
"4", {\ar"5"},
"4", {\ar"7"},
"8", {\ar"4"},
"5", {\ar"8"},
"6", {\ar"7"},
"7", {\ar"8"},
\end{xy}
\]
It admits the maximal green sequence $3,6,9,2,5,8,1,4,7,3,6,9,2,5,8,3,6,9$.
In this case, the invariant $DT_Q$ is of infinite order.

% \acknowledgements
% \label{sec:ack}

\nocite{*}
\bibliographystyle{abbrvnat}
% use the following instead if you encounter problems 
%\bibliographystyle{alpha}
%\bibliography{Keller-FPSAC-2013-v2}

\begin{thebibliography}{60}
\providecommand{\natexlab}[1]{#1}
\providecommand{\url}[1]{\texttt{#1}}
\expandafter\ifx\csname urlstyle\endcsname\relax
  \providecommand{\doi}[1]{doi: #1}\else
  \providecommand{\doi}{doi: \begingroup \urlstyle{rm}\Url}\fi

\bibitem[Adachi et~al.(2012)Adachi, Iyama, and Reiten]{AdachiIyamaReiten12}
T.~Adachi, O.~Iyama, and I.~Reiten.
\newblock $ \tau$-tilting theory.
\newblock arXiv:1210.1036 [math.RT], 2012.

\bibitem[Alim et~al.(2011{\natexlab{a}})Alim, Cecotti, C\'ordova, Espahbodi,
  Rastogi, and Vafa]{CecottiEtAl11}
M.~Alim, S.~Cecotti, C.~C\'ordova, S.~Espahbodi, A.~Rastogi, and C.~Vafa.
\newblock $\mathcal{N}=2$ quantum field theories and their {BPS} quivers.
\newblock arXiv:1112.3984 [hep-th], 2011{\natexlab{a}}.

\bibitem[Alim et~al.(2011{\natexlab{b}})Alim, Cecotti, C\'ordova, Espahbodi,
  Rastogi, and Vafa]{CecottiEtAl11a}
M.~Alim, S.~Cecotti, C.~C\'ordova, S.~Espahbodi, A.~Rastogi, and C.~Vafa.
\newblock {BPS} quivers and spectra of complete $\mathcal{N}=2$ quantum field
  theories.
\newblock arXiv:1109.4941 [hep-th], 2011{\natexlab{b}}.

\bibitem[Berenstein and Zelevinsky(2005)]{BerensteinZelevinsky05}
A.~Berenstein and A.~Zelevinsky.
\newblock Quantum cluster algebras.
\newblock \emph{Adv. Math.}, 195\penalty0 (2):\penalty0 405--455, 2005.

\bibitem[Berenstein et~al.(2005)Berenstein, Fomin, and
  Zelevinsky]{BerensteinFominZelevinsky05}
A.~Berenstein, S.~Fomin, and A.~Zelevinsky.
\newblock Cluster algebras. {III}. {U}pper bounds and double {B}ruhat cells.
\newblock \emph{Duke Math. J.}, 126\penalty0 (1):\penalty0 1--52, 2005.

\bibitem[Br\"ustle et~al.(2012)Br\"ustle, Dupont, and
  P\'erotin]{BruestleDupontPerotin12}
T.~Br\"ustle, G.~Dupont, and M.~P\'erotin.
\newblock On maximal green sequences.
\newblock arXiv:1205.2050 [math.RT], 2012.

\bibitem[Caldero and Chapoton(2006)]{CalderoChapoton06}
P.~Caldero and F.~Chapoton.
\newblock Cluster algebras as {H}all algebras of quiver representations.
\newblock \emph{Comment. Math. Helv.}, 81\penalty0 (3):\penalty0 595--616,
  2006.

\bibitem[Cecotti and Vafa(2011)]{CecottiVafa11}
S.~Cecotti and C.~Vafa.
\newblock Classification of complete {$N=2$} supersymmetric field theories in
  $4$ dimensions.
\newblock arXiv: 1103.5832 [hep-th], 2011.

\bibitem[Cecotti et~al.(2010)Cecotti, Neitzke, and Vafa]{CecottiNeitzkeVafa10}
S.~Cecotti, A.~Neitzke, and C.~Vafa.
\newblock {$R$-twisting and $4d$/$2d$-Correspondences}.
\newblock arXiv:10063435 [physics.hep-th], 2010.

\bibitem[Cecotti et~al.(2011)Cecotti, C\'{o}rdova, and
  Vafa]{CecottiCordovaVafa11}
S.~Cecotti, C.~C\'{o}rdova, and C.~Vafa.
\newblock Braids, walls and mirrors.
\newblock arXiv:1110.2115 [hep-th], 2011.

\bibitem[Cerulli~Irelli et~al.(2012)Cerulli~Irelli, Keller, Labardini-Fragoso,
  and Plamondon]{CerulliKellerLabardiniPlamondon12}
G.~Cerulli~Irelli, B.~Keller, D.~Labardini-Fragoso, and P.-G. Plamondon.
\newblock Linear independence of cluster monomials for skew-symmetric cluster
  algebras.
\newblock arXiv:1203.1307 [math.RT], to appear in Compositio, 2012.

\bibitem[Chapoton et~al.(2002)Chapoton, Fomin, and
  Zelevinsky]{ChapotonFominZelevinsky02}
F.~Chapoton, S.~Fomin, and A.~Zelevinsky.
\newblock Polytopal realizations of generalized associahedra.
\newblock \emph{Canad. Math. Bull.}, 45\penalty0 (4):\penalty0 537--566, 2002.
\newblock Dedicated to Robert V.\ Moody.

\bibitem[Derksen et~al.(2008)Derksen, Weyman, and
  Zelevinsky]{DerksenWeymanZelevinsky08}
H.~Derksen, J.~Weyman, and A.~Zelevinsky.
\newblock Quivers with potentials and their representations {I}: {Mutations}.
\newblock \emph{Selecta Mathematica}, 14:\penalty0 59--119, 2008.

\bibitem[Derksen et~al.(2010)Derksen, Weyman, and
  Zelevinsky]{DerksenWeymanZelevinsky10}
H.~Derksen, J.~Weyman, and A.~Zelevinsky.
\newblock Quivers with potentials and their representations {II}: {Applications
  to cluster algebras}.
\newblock \emph{J.~Amer.~Math.~Soc.}, 23:\penalty0 749--790, 2010.

\bibitem[Faddeev and Volkov(1993)]{FaddeevVolkov93}
L.~Faddeev and A.~Y. Volkov.
\newblock Abelian current algebra and the {V}irasoro algebra on the lattice.
\newblock \emph{Phys. Lett. B}, 315\penalty0 (3-4):\penalty0 311--318, 1993.

\bibitem[Faddeev and Kashaev(1994)]{FaddeevKashaev94}
L.~D. Faddeev and R.~M. Kashaev.
\newblock Quantum dilogarithm.
\newblock \emph{Modern Phys. Lett. A}, 9\penalty0 (5):\penalty0 427--434, 1994.

\bibitem[Felikson et~al.(2012{\natexlab{a}})Felikson, Shapiro, and
  Tumarkin]{FeliksonShapiroTumarkin12}
A.~Felikson, M.~Shapiro, and P.~Tumarkin.
\newblock Skew-symmetric cluster algebras of finite mutation type.
\newblock \emph{J. Eur. Math. Soc. (JEMS)}, 14\penalty0 (4):\penalty0
  1135--1180, 2012{\natexlab{a}}.
\newblock ISSN 1435-9855.
\newblock \doi{10.4171/JEMS/329}.
\newblock URL \url{http://dx.doi.org/10.4171/JEMS/329}.

\bibitem[Felikson et~al.(2012{\natexlab{b}})Felikson, Shapiro, and
  Tumarkin]{FeliksonShapiroTumarkin12a}
A.~Felikson, M.~Shapiro, and P.~Tumarkin.
\newblock Cluster algebras of finite mutation type via unfoldings.
\newblock \emph{Int. Math. Res. Not. IMRN}, \penalty0 (8):\penalty0 1768--1804,
  2012{\natexlab{b}}.
\newblock ISSN 1073-7928.
\newblock \doi{10.1093/imrn/rnr072}.
\newblock URL \url{http://dx.doi.org/10.1093/imrn/rnr072}.

\bibitem[Fock and Goncharov(2009{\natexlab{a}})]{FockGoncharov09}
V.~V. Fock and A.~B. Goncharov.
\newblock Cluster ensembles, quantization and the dilogarithm.
\newblock \emph{Annales scientifiques de l'ENS}, 42\penalty0 (6):\penalty0
  865--930, 2009{\natexlab{a}}.

\bibitem[Fock and Goncharov(2009{\natexlab{b}})]{FockGoncharov09a}
V.~V. Fock and A.~B. Goncharov.
\newblock The quantum dilogarithm and representations of quantum cluster
  varieties.
\newblock \emph{Invent. Math.}, 175\penalty0 (2):\penalty0 223--286,
  2009{\natexlab{b}}.

\bibitem[Fomin(2002)]{Fomin07}
S.~Fomin.
\newblock Cluster algebras portal.
\newblock \verb"www.math.lsa.umich.edu/~fomin/"\verb"cluster.html", 2002.

\bibitem[Fomin(2010)]{Fomin10}
S.~Fomin.
\newblock Total positivity and cluster algebras.
\newblock In \emph{Proceedings of the {I}nternational {C}ongress of
  {M}athematicians. {V}olume {II}}, pages 125--145, New Delhi, 2010. Hindustan
  Book Agency.

\bibitem[Fomin and Zelevinsky(2002)]{FominZelevinsky02}
S.~Fomin and A.~Zelevinsky.
\newblock Cluster algebras. {I}. {F}oundations.
\newblock \emph{J. Amer. Math. Soc.}, 15\penalty0 (2):\penalty0 497--529
  (electronic), 2002.

\bibitem[Fomin and Zelevinsky(2003)]{FominZelevinsky03b}
S.~Fomin and A.~Zelevinsky.
\newblock {$Y$}-systems and generalized associahedra.
\newblock \emph{Ann. of Math. (2)}, 158\penalty0 (3):\penalty0 977--1018, 2003.

\bibitem[Fomin and Zelevinsky(2007)]{FominZelevinsky07}
S.~Fomin and A.~Zelevinsky.
\newblock Cluster algebras {IV}: {C}oefficients.
\newblock \emph{Compositio Mathematica}, 143:\penalty0 112--164, 2007.

\bibitem[Fomin et~al.(2008)Fomin, Shapiro, and
  Thurston]{FominShapiroThurston08}
S.~Fomin, M.~Shapiro, and D.~Thurston.
\newblock Cluster algebras and triangulated surfaces. {I}. {C}luster complexes.
\newblock \emph{Acta Math.}, 201\penalty0 (1):\penalty0 83--146, 2008.

\bibitem[Gaiotto et~al.(2008)Gaiotto, Moore, and
  Neitzke]{GaiottoMooreNeitzke08}
D.~Gaiotto, G.~W. Moore, and A.~Neitzke.
\newblock Four-dimensional wall-crossing via three-dimensional field theory.
\newblock arXiv:0807.4723 [physics.hep-th], 2008.

\bibitem[Gaiotto et~al.(2009)Gaiotto, Moore, and
  Neitzke]{GaiottoMooreNeitzke09}
D.~Gaiotto, G.~W. Moore, and A.~Neitzke.
\newblock Wall-crossing, {H}itchin systems, and the {WKB} approximation.
\newblock arXiv:0907.3987 [physics.hep-th], 2009.

\bibitem[Gaiotto et~al.(2010{\natexlab{a}})Gaiotto, Moore, and
  Neitzke]{GaiottoMooreNeitzke10}
D.~Gaiotto, G.~W. Moore, and A.~Neitzke.
\newblock {Framed BPS States}.
\newblock arXiv:1006.0146 [physics.hep-th], 2010{\natexlab{a}}.

\bibitem[Gaiotto et~al.(2010{\natexlab{b}})Gaiotto, Moore, and
  Neitzke]{GaiottoMooreNeitzke10a}
D.~Gaiotto, G.~W. Moore, and A.~Neitzke.
\newblock Four-dimensional wall-crossing via three-dimensional field theory.
\newblock \emph{Comm. Math. Phys.}, 299\penalty0 (1):\penalty0 163--224,
  2010{\natexlab{b}}.

\bibitem[Gei\ss et~al.(2011)Gei\ss, Leclerc, and
  Schr{\"o}er]{GeissLeclercSchroeer11b}
C.~Gei\ss, B.~Leclerc, and J.~Schr{\"o}er.
\newblock Kac-{M}oody groups and cluster algebras.
\newblock \emph{Adv. Math.}, 228\penalty0 (1):\penalty0 329--433, 2011.
\newblock ISSN 0001-8708.
\newblock \doi{10.1016/j.aim.2011.05.011}.
\newblock URL \url{http://dx.doi.org/10.1016/j.aim.2011.05.011}.

\bibitem[Joyce and Song(2009)]{JoyceSong09}
D.~Joyce and Y.~Song.
\newblock A theory of generalized {Donaldson-Thomas} invariants. {II}.
  {M}ultiplicative identities for {B}ehrend functions.
\newblock arXiv:0901.2872 [math.AG], 2009.

\bibitem[Keller(2006)]{KellerQuiverMutationApplet}
B.~Keller.
\newblock Quiver mutation in {J}ava.
\newblock Java applet available at the author's home page, 2006.

\bibitem[Keller(2010)]{Keller10b}
B.~Keller.
\newblock Cluster algebras, quiver representations and triangulated categories.
\newblock In T.~Holm, P.~J{\o}rgensen, and R.~Rouquier, editors,
  \emph{Triangulated categories}, volume 375 of \emph{London Mathematical
  Society Lecture Note Series}, pages 76--160. Cambridge University Press,
  2010.

\bibitem[Keller(2011)]{Keller11c}
B.~Keller.
\newblock On cluster theory and quantum dilogarithm identities.
\newblock In \emph{Representations of algebras and related topics}, EMS Ser.
  Congr. Rep., pages 85--116. Eur. Math. Soc., Z\"urich, 2011.
\newblock \doi{10.4171/101-1/3}.
\newblock URL \url{http://dx.doi.org/10.4171/101-1/3}.

\bibitem[Keller(2012)]{Keller12}
B.~Keller.
\newblock Cluster algebras and derived categories.
\newblock arXiv:1202.4161 [math.RT], 2012.

\bibitem[King and Qiu(2011)]{KingQiu11}
A.~King and Y.~Qiu.
\newblock Exchange graphs of acyclic {C}alabi--{Y}au categories.
\newblock arXiv:1109.2924 [math.RT], 2011.

\bibitem[Koenig and Yang(2012)]{KoenigYang12}
S.~Koenig and D.~Yang.
\newblock Silting objects, simple-minded collections, $t$-structures and
  co-$t$-structures for finite-dimensional algebras.
\newblock arXiv:1203.5657 [math.RT], 2012.

\bibitem[Kontsevich and Soibelman(2008)]{KontsevichSoibelman08}
M.~Kontsevich and Y.~Soibelman.
\newblock Stability structures, {D}onaldson-{T}homas invariants and cluster
  transformations.
\newblock arXiv:0811.2435 [math.AG], 2008.

\bibitem[Kontsevich and Soibelman(2010)]{KontsevichSoibelman10}
M.~Kontsevich and Y.~Soibelman.
\newblock Cohomological {H}all algebra, exponential {H}odge structures and
  motivic {D}onaldson-{T}homas invariants.
\newblock arXiv:1006.2706 [math.AG], 2010.

\bibitem[Ladkani(2007)]{Ladkani07}
S.~Ladkani.
\newblock Universal derived equivalences of posets of tilting modules.
\newblock arXiv:0708.1287 [math.RT], 2007.

\bibitem[Leclerc(2010)]{Leclerc10}
B.~Leclerc.
\newblock Cluster algebras and representation theory.
\newblock In \emph{Proceedings of the {I}nternational {C}ongress of
  {M}athematicians. {V}olume {IV}}, pages 2471--2488, New Delhi, 2010.
  Hindustan Book Agency.

\bibitem[Musiker and Stump(2011)]{MusikerStump11}
G.~Musiker and C.~Stump.
\newblock A compendium on the cluster algebra and quiver package in {S}age.
\newblock arXiv:1102.4844 [math.CO], 2011.

\bibitem[Nagao(2010)]{Nagao10}
K.~Nagao.
\newblock {Donaldson-Thomas theory and cluster algebras}.
\newblock arXiv:1002.4884 [math.AG], 2010.

\bibitem[Nagao(2011)]{Nagao11}
K.~Nagao.
\newblock Quantum dilogarithm identities.
\newblock \emph{RIMS Kokyuroku Bessatsu}, B28:\penalty0 164--170, 2011.

\bibitem[Nakanishi(2012)]{Nakanishi12}
T.~Nakanishi.
\newblock Tropicalization method in cluster algebras.
\newblock \emph{Contemp. Math.}, 580:\penalty0 95--115, 2012.

\bibitem[Plamondon(2011)]{Plamondon11a}
P.-G. Plamondon.
\newblock Cluster algebras via cluster categories with infinite-dimensional
  morphism spaces.
\newblock \emph{Compositio Mathematica}, 147:\penalty0 1921--1954, 2011.

\bibitem[Qiu(2011)]{Qiu11}
Y.~Qiu.
\newblock Stability conditions and quantum dilogarithm identities for {D}ynkin
  quivers.
\newblock arXiv:1111.1010 [math.AG], 2011.

\bibitem[Reading(2006)]{Reading06}
N.~Reading.
\newblock Cambrian lattices.
\newblock \emph{Adv. Math.}, 205\penalty0 (2):\penalty0 313--353, 2006.
\newblock ISSN 0001-8708.
\newblock \doi{10.1016/j.aim.2005.07.010}.
\newblock URL \url{http://dx.doi.org/10.1016/j.aim.2005.07.010}.

\bibitem[Reineke(2010)]{Reineke10}
M.~Reineke.
\newblock Poisson automorphisms and quiver moduli.
\newblock \emph{J. Inst. Math. Jussieu}, 9\penalty0 (3):\penalty0 653--667,
  2010.

\bibitem[Reineke(2011)]{Reineke11}
M.~Reineke.
\newblock Cohomology of quiver moduli, functional equations, and integrality of
  {D}onaldson-{T}homas type invariants.
\newblock \emph{Compos. Math.}, 147\penalty0 (3):\penalty0 943--964, 2011.
\newblock ISSN 0010-437X.
\newblock \doi{10.1112/S0010437X1000521X}.
\newblock URL \url{http://dx.doi.org/10.1112/S0010437X1000521X}.

\bibitem[Reiten(2010)]{Reiten10a}
I.~Reiten.
\newblock Cluster categories.
\newblock In \emph{Proceedings of the {I}nternational {C}ongress of
  {M}athematicians. {V}olume {I}}, pages 558--594, New Delhi, 2010. Hindustan
  Book Agency.

\bibitem[Seiberg(1995)]{Seiberg95}
N.~Seiberg.
\newblock Electric-magnetic duality in supersymmetric non-abelian gauge
  theories.
\newblock \emph{Nuclear Phys. B}, 435\penalty0 (1-2):\penalty0 129--146, 1995.
\newblock arXiv:hep-th/9411149.

\bibitem[Szendr{\H{o}}i(2008)]{Szendroi08}
B.~Szendr{\H{o}}i.
\newblock Non-commutative {D}onaldson-{T}homas invariants and the conifold.
\newblock \emph{Geom. Topol.}, 12\penalty0 (2):\penalty0 1171--1202, 2008.
\newblock \doi{10.2140/gt.2008.12.1171}.
\newblock URL \url{http://dx.doi.org/10.2140/gt.2008.12.1171}.

\bibitem[Thomas(2011)]{Thomas11}
H.~Thomas.
\newblock The {T}amari lattice as it arises in quiver representations.
\newblock arXiv:1110.3040 [math.RT], 2011.

\bibitem[Volkov(2012)]{Volkov12}
A.~Y. Volkov.
\newblock Pentagon identity revisited.
\newblock \emph{Int. Math. Res. Not. IMRN}, \penalty0 (20):\penalty0
  4619--4624, 2012.
\newblock ISSN 1073-7928.

\bibitem[Williams(2012)]{Williams12}
L.~Williams.
\newblock Cluster algebras: an introduction.
\newblock arXiv:1212.6263 [math.RA], 2012.

\bibitem[Xie(2012)]{Xie12}
D.~Xie.
\newblock {BPS} spectrum, wall crossing and quantum dilogarithm identity.
\newblock arXiv:1211.7071 [physics.hep-th], 2012.

\bibitem[Zagier(1991)]{Zagier91}
D.~Zagier.
\newblock Polylogarithms, {D}edekind zeta functions and the algebraic
  {$K$}-theory of fields.
\newblock In \emph{Arithmetic algebraic geometry ({T}exel, 1989)}, volume~89 of
  \emph{Progr. Math.}, pages 391--430. Birkh{\"a}user Boston, Boston, MA, 1991.

\bibitem[Zamolodchikov(1991)]{Zamolodchikov91}
A.~B. Zamolodchikov.
\newblock On the thermodynamic {B}ethe ansatz equations for reflectionless
  {$ADE$} scattering theories.
\newblock \emph{Phys. Lett. B}, 253\penalty0 (3-4):\penalty0 391--394, 1991.

\end{thebibliography}
\label{sec:biblio}

\end{document}